\documentclass[11pt]{amsart}
\usepackage{amssymb}
\usepackage{amscd}
\newtheorem{theo}{Theorem}

\newtheorem{prop}{Proposition}
\newtheorem{rem}{Remark}
\newtheorem{lem}{Lemma}

\newenvironment{pf}{{\bf Proof. }}{\hfill $\square$\medskip}

\newcommand{\id}{\operatorname{id}}
\newcommand{\C}{{\mathbb C}}

\newcommand{\CP}{{\mathbb C \mathbb P}}

\newcommand{\e}{\operatorname{e}}
\newcommand{\I}{\operatorname{i}}
\newcommand{\Aut}{\operatorname{Aut}}

\newcommand{\SL}{\operatorname{SL}}
\newcommand{\PSL}{\operatorname{PSL}}

\newcommand{\pd}[2][]{\tfrac{\partial^{#1}}{\partial #2}}

\def\sideremark#1{\ifvmode\leavevmode\fi\vadjust{
\vbox to0pt{\hbox to 0pt{\hskip\hsize\hskip1em
\vbox{\hsize3cm\tiny\raggedright\pretolerance10000
\noindent #1\hfill}\hss}\vbox to8pt{\vfil}\vss}}}
\begin{document}
\title[Elliptic CR-manifolds with additional symmetries]{Elliptic CR-manifolds and shear invariant
ODE with additional symmetries}
\author{Vladimir Ezhov and Gerd Schmalz}
\date{}
\address{(Ezhov) School of Mathematics and Statistics,
University of South Australia, Mawson Lakes Boulevard, Mawson
Lakes, 5095 South Australia} \email{vladimir.ejov@unisa.edu.au}

\address{(Schmalz) School of Mathematics, Statistics and Computer Science\\
University of New England\\Armidale, NSW
2351\\Australia//Uniwersytet Warmi\'nsko-Mazurski, Wydzia{\l}
Matematyki-Informatyki, ul. \.Zo{\l}nierska 14, 10-561 Olsztyn,
Poland} \email{gerd@turing.une.edu.au}
\begin{abstract}
We classify the ODEs that correspond to elliptic CR-mani\-folds with
maximal isotropy. It follows that the dimension of the isotropy
group of an elliptic CR-manifold can be only 10 (for the
quadric), 4 (for the listed examples) or less. This is in
contrast with the situation of hyperbolic CR-manifolds, where the
dimension can be 10 (for the quadric), 6 or 5 (for
semi-quadrics) or less than 4. We also prove that, for all
elliptic CR-manifolds with non-linearizable istropy group, except
for two special manifolds, the points with non-linearizable
isotropy form exactly some complex curve on the manifold.
\end{abstract}

\maketitle

\section{Introduction}

In \cite{ES03} the authors used a correspondence between so-called
torsion-free elliptic CR-manifolds and complex second order ODE to
describe elliptic CR-manifolds with non-linearizable isotropy.
This description was based on an investigation of ODE's with a
shear symmetry $y \frac{\partial}{\partial x}$ on the $x,y$-plane
near the singularity $(0,0)$.

The major aim of this paper is to describe elliptic CR-manifolds
with big isotropy. We will show that the maximal dimension of the
isotropy for non-quadratic elliptic CR-manifolds is $4$ and is
attained exactly for manifolds that correspond to the ODE's
\begin{align*}
y''&=y^k (y-xy')^3\\
y''&=y^{\ell} y'(y-xy')^2+C y^{2\ell +2}(y-xy')^3
\end{align*}

where $k,\ell$ are non-negative integers and $C$ is a complex
constant. Thus, according to earlier results by the authors
\cite{ESs},  the possible dimensions of the isotropy of elliptic
CR-manifolds are $10,4,3,2,1,0$. This is somewhat unexpected,
because the corresponding numbers for analogous hyperbolic
manifolds are $10,6,5,3,2,1,0$ (see \cite{Hab}).

In Section \ref{sec4} we represent open parts of elliptic
CR-manifolds with non-linearizable isotropy as copies of
$\SL(2,\C)$ with the standard action of subgroups of $\SL(2,\C)$.

In Section \ref{secdu} we show that the duality of ODE that
results from switching the r\^oles of variables $x,y$ and
parameters $c_1,c_2$ corresponds simply to switching to the
complex conjugate CR-manifold. We demonstrate this feature for the
exceptional quartic.

Section \ref{onenon} is devoted to shear-invariant elliptic
CR-manifolds with one additional non-isotropic symmetry. We show
that these manifolds coincide with the manifolds obtained in
Section \ref{oneiso} for a different choice of the reference
point.

In Section \ref{twosy} we conclude that the quartic is the only
shear invariant elliptic CR-manifold with 6-dimensional
automorphism group.

Finally, in Section \ref{nonlin} we show that the quadric and the
quartic are characterized by the property that the points with
non-linearizable isotropy fill more than a complex curve, whereas
in all other cases, they fill exactly a complex curve.

\section{Preliminaries}

Let $M$ be a CR-manifold $M$ of CR-dimension two and
CR-codimension two, i.e., $M$ is a $6$-dimensional manifold with a
$4$-dimensional distribution $D \subset TM$ and a smooth field of
endomorphisms $J_x: D_x \to D_x$ with $J_x^2=-\id$. The Levi form
at $x\in M$ is a bilinear mapping
$${\mathcal L}_x: D_x \times D_x \to T_xM/ D_x.$$

${\mathcal L}_x(X,Y)$ is defined as the bracket of two sections
$\tilde{X}, \tilde{Y}$ of $D$ that extend $X,Y$, followed by the
natural projection $\pi: T_x \to T_x / D_x $.

$M$ is called elliptic if any {\em real} linear combination of the two
scalar components of $\mathcal L$ is a non-degenerate bilinear
form. It follows that there exist two mutually conjugate {\em complex}
degenerate combinations. Its null vectors define a canonical
splitting  $D_x= D^+_x\oplus D^-_x$. For a pair of sections
$\tilde{X}$ in $D^+$ and $\tilde{Y}$ in $D^-$ the vectors
$(\tilde{X}_x, \bar{\tilde{Y}}_x, [\tilde{X},\bar{\tilde{Y}}]_x)$
define a complex structure on $T_xM$.

We assume that ${\mathcal L}_x (J_x X, J_x Y)={\mathcal L}_x
(X,Y)$ for all $x\in M$, i.e., $M$ is partially integrable.

For partially integrable elliptic CR-manifolds a Cartan connection
was constructed in \cite{CS} (see also \cite{SS,SSp}). In this
paper we consider only elliptic manifolds, whose so-called
torsion-part of the Cartan curvature vanishes. This algebraic
condition is equivalent to the following geometric properties:

\begin{enumerate}
\item $M$ is embeddable,
\item the line bundles $D^+$ and $D^-$ are integrable,
\item the canonical almost complex structure is integrable.
\end{enumerate}

It follows that $M$ must be real analytic. For a smooth embedded
elliptic CR-manifold vanishing of the torsion at a point $x\in M$
can also be expressed by the equivalent condition that $M$ has
contact of third order with its osculating quadric at $x$ (see
\cite{SSa}).

We have

\begin{prop}\label{prop1}
There is a 1-1 correspondence between
\begin{enumerate}
\item Torsionfree elliptic CR-manifolds of CR-dimension two and
CR-codimension two,
\item Complex 3-folds with two holomorphic direction fields that span a non-involutive distribution,
\item Complex second order ODE.
\end{enumerate}

\end{prop}

\begin{pf}
If $\tilde{M}$ is a complex $3$-fold with a pair of non-involutive
direction fields one can introduce local coordinates $x,y,p$ such
that $\mathfrak Z_1=\pd{p}{}$ and $\mathfrak Z_2=\pd{x}{} +p
\pd{y}{} +B(x,y,p)\pd{p}{}$ (see \cite{ES03}). This allows
an interpretation of a local part of $\tilde{M}$ as a chart of the
projectivized tangent bundle over $\C^2$ with coordinates $x,y$ in
the base and $p=\frac{dy}{dx}$ in the fibre. The projections
of the integral curves of $\mathfrak Z_2$ are then nothing but the
integral curves of $y''=B(x,y,y')$ in $\C^2$. Vice versa, the
lifts of integral curves of a second order ODE to the
projectivized tangent bundle define a direction field $\mathfrak
Z_2$ that does not commute with $\mathfrak Z_1=\pd{p}{}$.

If $M$ is a torsionfree elliptic CR-manifold then $M$ has an
integrable almost complex structure and two holomorphic direction
fields that generate $D^+, \bar{D}^-$. Vice versa, if $\tilde{M}$
is a complex $3$-fold with holomorphic direction fields $\mathfrak
Z_1, \mathfrak Z_2$ then $D_x$ can be defined as the span of these
direction fields. $J_x$ is defined by $J_x \mathfrak Z_{1,x}= \I
\mathfrak Z_{1,x}$ and $J_x \mathfrak Z_{2,x}= -\I \mathfrak Z_{2,x}$.
Non-involutivity of the two direction fields is equivalent to the
ellipticity of the Levi form.

It remains to show that the obtained CR-manifold $M$ is
torsionfree. It is convenient to represent $M$ as an embedded
CR-submanifold of $\C^4$. We look for four independent coordinate
functions that are annihilated by
$$\bar{\mathfrak Z}_1=\pd{\bar{p}}{}, \quad \mathfrak Z_2=\pd{x}{}
+p \pd{y}{} +B(x,y,p)\pd{p}{}$$

Two obvious solutions are $z_2=\bar{x}$ and $w_2=\bar{y}$. We need
two additional coordinate functions of the form $f(x,y,p)$. Thus,
we have to solve
$$\pd{x}{f}+p \pd{y}{f} +B(x,y,p)\pd{p}{f}=0.$$

The characteristic equation of this PDE is
$$\dot{x}=1, \quad \dot{y}=p, \quad \dot{p}=B(x,y,p).$$

It is equivalent to $\ddot{y}=B(t,y,\dot{y})$. Let
$$ x=t,\quad y=\phi(t,C_1,C_2), \quad p=\dot{\phi}(t,C_1,C_2)$$
be the characteristic curves. Then the desired coordinate functions are
$z_1=C_1(x,y,p)$ and $z_2=C_2(x,y,p)$. In $\C^4$ with coordinates
$z_1,z_2,w_1,w_2$ the equation of the manifold $M$ takes the form
$$\bar{w}_2=\phi(\bar{z}_2,z_1,w_1).$$
$M$ has two foliations: into holomorphic curves (for $\bar{z}_2,
\bar{w}_2$ fixed) and into antiholomorphic curves (for $z_1, w_1$
fixed). The tangent spaces to the curves that pass through a given
point span the maximal complex subspace of the tangent space of
$M$ at this point. The corresponding directions annihilate
degenerate complex linear combination of the components of the
Levi form. Thus, they provide the canonical splitting. By
construction, the corresponding line bundles are integrable. The
induced almost complex structure is the one that is obtained by
adopting $z_1,\bar{z}_2,w_1,\bar{w}_2$ as holomorphic coordinates
in the ambient space. Therefore, it is clearly integrable.
\end{pf}

\begin{rem}
The embedding constructed in the proof of Proposition \ref{prop1}
has the property that the two canonical foliations coincide with
the foliations into the fibres of the projections to the
$z_1,w_1$-plane and the $z_2,w_2$-plane, respectively.
\end{rem}

It was proved in \cite{ES03} that elliptic CR-manifolds with
non-linearizable isotropy group are in 1-1 correspondence with
shear invariant second order ODE. Such ODE can be represented by
\begin{equation}\label{ode}
y''=B(x,y,y')=f_0(y) (y-xy')^3 +f_1(y) y' (y-xy')^2,
\end{equation}

where two ODE are equivalent if and only if there is a mapping
$$(x,y) \to \left(\frac{c_1x}{1-cy}, \frac{c_2y}{1-cy}\right)$$

that takes one to the other. A finer classification can be
obtained if we take into account possible additional symmetries.
Sophus Lie \cite{Lie83} classified second order ODE with one-,
two-, and three-dimensional symmetry groups. The difference of our
approach is that we are interested in fixed points of the
automorphisms, whereas Lie always chooses a point, where one of
the symmetries is a translation $\pd{y}{}$. In our situation one
of the symmetries is the shear $y\pd{x}{}$. Our choice of the
canonical symmetries and regularity of the ODE at the reference
point imply that $B$ is a third order polynomial with respect to
$y'$ and $x$.

\section{Classification of shear invariant ODE with 4-dimensional
isotropy}\label{oneiso}

If there is only one (up to scale) shear symmetry of a shear
invariant ODE then it can be used as an invariant. On the other
hand, as it is known from \cite{ES03}, all ODEs with more than one shear
can be written as
$$y''=\frac{K(y-xy')^3}{(1-cy)^3}.$$
If we exclude these then any additional isotropic symmetry of the
ODE $y''=B(x,y,y')$ must preserve the single shear symmetry and,
consequently, must have the form
\begin{equation}\label{iaut}
((\phi(y)+a)x +\psi(y))\pd{x}{} + \phi(y)y \pd{y}{}.
\end{equation}

The general equation for infinitesimal symmetries $\xi\pd{x}{}
+\eta \pd{y}{}$ is
\begin{multline}\label{inf}
\xi \frac{\partial B}{\partial x} + \eta  \frac{\partial
B}{\partial y} +\phi  \frac{\partial B}{\partial p} +\left(2
\frac{\partial \xi}{\partial
    x} + 3p  \frac{\partial \xi}{\partial y} - \frac{\partial \eta}{\partial
    y}\right) B -\\ -\frac{\partial^2 \eta}{(\partial x)^2} + p
\left(\frac{\partial^2 \xi}{(\partial x)^2} - 2 \frac{\partial^2
\eta}{\partial
    x \partial y}  \right) + p^2 \left( 2 \frac{\partial^2 \xi}{\partial
    x \partial y}- \frac{\partial^2 \eta}{(\partial y)^2}\right) + p^3
\frac{\partial^2 \xi}{(\partial y)^2}=0.
\end{multline}

We plug in $B$ from (\ref{ode}) and the infinitesimal automorphism
(\ref{iaut}). The component of degree $3$ in $p$ and degree $0$ in
$x$ in equation (\ref{inf}) immediately implies $\psi''=0$. Since
we are here interested only in isotropic automorphisms and since
we know that the shear $y\pd{x}{}$ is an automorphism we may
assume $\psi=0$.

The component of degree $3$ in $p$ and degree $1$ in $x$ in
equation (\ref{inf}) yields now $\phi''=0$, thus $\phi=\beta_1 +
\alpha_3 y$

From the components of degree $3$ in $p$ and degree $2$ and $3$ in
$x$ we get
\begin{align*}
af_1+3\beta_1 f_1+3\alpha_3 y f_1+ \beta_1 y
f_1'+\alpha_3y^2f_1'&=0\\
2af_0 +4\beta_1 f_0+ 3y\alpha_3f_0+ y\beta_1 f_0'+ y^2\alpha_3
f_0'&=0.
\end{align*}

If $f_0=\sum_{n=k}^\infty b_n y^n$ and $f_1=\sum_{n=\ell}^\infty
c_n y^n$ then
\begin{align*}
(a+(n+3)\beta_1) c_n + (n+2)\alpha_3 c_{n-1}&=0\\
(2a+(n+4)\beta_1) b_n + (n+2)\alpha_3 b_{n-1}&=0.
\end{align*}

The first equation for $n=\ell$ and the second equation for $n=k$
give rise to a linear system that implies $\beta_1=a=0$ and,
consequently, $\alpha_3=0$, unless $k=2\ell+2$, or either $f_0=0$
or $f_1=0$.

From the recursive formulae we find
\begin{align*}
f_0&= C_1 \,y^k (1-cy)^{-k-3}\\
f_1&= C_2 \,y^\ell (1-cy)^{-\ell-3}
\end{align*}

By applying a transformation $x_1=\frac{c_1 x}{1-cy}$,
$y_1=\frac{c_2 y}{1-cy}$ this can be reduced to one of the
following two series of ODE
\begin{align}\label{vier1}
y''&=y^k (y-xy')^3\\\label{vier3} y''&=y^{\ell} y'(y-xy')^2+C
y^{2\ell +2} (y-xy')^3
\end{align}

where $k,\ell$ are non-negative integers and $C$ is a complex
constant. According to Theorem 3 in \cite{ES03} these ODE are
pairwise non-equivalent.

The additional symmetry is
$$(k+2) x \pd{x}{} -2 y \pd{y}{}, \quad \text{resp.} \quad (\ell+2) x \pd{x}{} -y \pd{y}{}$$

The corresponding CR-manifolds are exactly the CR-manifolds with
an isotropy group of real dimension $4$.

We conclude
\begin{theo}
The isotropy group of an elliptic CR-manifold has
\begin{enumerate}
\item dimension $10$ if and only if it is equivalent to the
quadric \item dimension $4$ if and only if it corresponds to one
of the ODE (\ref{vier1}) - (\ref{vier3}) \item dimension $\le 3$
in all other cases.
\end{enumerate}
\end{theo}

\begin{pf}
Statements (1) and (2) follow from the obtained classification.
Statement (3) was proved in \cite{ESs}.
\end{pf}

\section{$\SL(2,\C)$ representation of the shear invariant
manifolds} \label{sec4}

Any shear invariant ODE
$$y''=f_0(y)(y-xy')^3 + f_1(y)y'(y-xy')^2$$

obviously admits the solutions $y=cx$ for any constant $c$. Thus,
the solution passing through $(x_0,y_0)\in \C^2_*= \C^2\setminus
\{(0,0)\}$ with slope $p_0=y_0/x_0$ is $y=p_0 x$.

Notice that the equation $y=p x$ describes a canonical section in
the trivial fibre bundle $\C^2_*\times \CP^1$, which is induced by
the tautological mapping
\begin{align*}
\tau:\, \C^2_*&\to \CP^1\\
(x,y) &\mapsto [x:y]
\end{align*}

By $M^*$ we denote the bundle with deleted section $\tau$.

Here we will give a representation of the part of the solution
manifold that corresponds to initial conditions $(x_0,y_0,p_0)\in
M^*$. $M^*$ can be identified with $\SL(2,\C)$ using the map
$$(x,y,p) \mapsto\begin{pmatrix}\frac{1}{y-xp} & x\\
\frac{p}{y-xp} &y \end{pmatrix} = \begin{pmatrix}\alpha & \beta\\
\gamma &\delta \end{pmatrix}.$$

The two distinguished direction fields now take the form
\begin{align*}
\mathfrak Z_1&= \beta \pd{\alpha}{} + \delta\pd{\gamma}{}\\
\mathfrak Z_2&=\alpha \pd{\beta}{} +  \gamma\pd{\delta}{} +
(f_0(\delta)+ \gamma f_1(\delta))\,\mathfrak Z_1
\end{align*}

The one-parametric action produced by the field $\mathfrak Z_1$ is
right multiplication with
$$\begin{pmatrix}1&0\\t&1\end{pmatrix}.$$

The second field generates a linear action only if $f_1\equiv 0$
and $f_0=const$, i.e., in the cases of a quadric ($f_0\equiv 0$)
or a quartic ($f_0\equiv 1$).

The shear symmetry is represented by
$$\theta=\gamma\pd{\alpha}{} + \delta\pd{\beta}{}$$

and produces left multiplication by
$$\begin{pmatrix}1&t\\0&1\end{pmatrix}.$$

In the quadric case $\mathfrak Z_2=\mathfrak Z_Q=\alpha
\pd{\beta}{} +  \gamma\pd{\delta}{} $ corresponds to right
multiplication with
$$\begin{pmatrix}1&t\\0&1\end{pmatrix}$$

and in the quartic case $\mathfrak Z_2=\mathfrak Z_Q+\mathfrak
Z_1$ to right multiplication with
$$\begin{pmatrix}\cosh t& \sinh t\\\sinh t&\cosh t\end{pmatrix}.$$

It is clear that in both cases these actions commute with the
complete left multiplication by $\SL(2,\C)$.

For manifolds with two isotropic symmetries the second (linear)
symmetry has the form
$$\mathfrak L=2\alpha\pd{\alpha}{} +(k+2) \beta\pd{\beta}{}
-(k+2) \gamma\pd{\gamma}{} -2 \delta\pd{\delta}{},$$

and respectively,
$$\mathfrak L=\alpha\pd{\alpha}{} +(\ell +2) \beta\pd{\beta}{}
-(\ell +2) \gamma\pd{\gamma}{} - \delta\pd{\delta}{},$$

It generates the one-parametric action
$$\begin{pmatrix}\alpha& \beta\\\gamma &\delta\end{pmatrix}\mapsto
\begin{pmatrix}t^{k+4}& 0\\0&t^{-k-4}\end{pmatrix}
\begin{pmatrix}\alpha& \beta\\\gamma &\delta\end{pmatrix}
\begin{pmatrix}t^{-k}& 0\\0&t^{k}\end{pmatrix},$$

and, respectively,
$$\begin{pmatrix}\alpha& \beta\\\gamma &\delta\end{pmatrix}\mapsto
\begin{pmatrix}t^{\ell+3}& 0\\0&t^{-\ell-3}\end{pmatrix}
\begin{pmatrix}\alpha& \beta\\\gamma &\delta\end{pmatrix}
\begin{pmatrix}t^{-\ell-1}& 0\\0&t^{\ell+1}\end{pmatrix},$$

Since $\mathfrak Z_1$ commutes with the left part of the action
and is mapped to $k \mathfrak Z_1$ (resp. $(\ell+1) \mathfrak
Z_1$) by the right part of the action we find
$$[\mathfrak L, \mathfrak Z_1]= k \mathfrak Z_1 \text{ resp. }  [\mathfrak L, \mathfrak Z_1]=
(\ell+1) \mathfrak Z_1.$$

For the second field
$$\mathfrak Z_2= \mathfrak Z_Q + F \mathfrak Z_1$$

we have
$$[\mathfrak L, \mathfrak Z_Q]= -k \mathfrak Z_Q \text{ resp. }  [\mathfrak L, \mathfrak Z_Q]=
(-\ell-1) \mathfrak Z_Q.$$

and
$$[\mathfrak L, F \mathfrak Z_1]= k F\mathfrak Z_1 + (\mathfrak L F) \mathfrak Z_1
\text{ resp. }  [\mathfrak L, F \mathfrak Z_1]= (\ell+1)
F\mathfrak Z_1 + (\mathfrak L F) \mathfrak Z_1.$$

In the first case this requires $\mathfrak L F= -2k F$, which is
satisfied for $F= \delta^k $. In the second case this requires
$\mathfrak L F= -2(\ell+1) F$ which is satisfied for combinations
of $\gamma\delta^\ell$ and $\delta^{2\ell+2}$.

\section{Dual ODE}\label{secdu}

A duality of ODE appears from the symmetry of interchanging the
distinguished direction fields $\mathfrak Z_1$ and $\mathfrak
Z_2$. This corresponds to interchanging the roles of the variables
$x,y$ and the parameters $c_1,c_2$ of the solutions. In terms of
the embedded CR-manifold this will be achieved by complex
conjugation. The symmetry group of the dual ODE clearly will be
isomorphic to the symmetry group of the initial ODE, though the
action is different. It follows that ODE corresponding to elliptic
CR-manifolds that are complexifications of real hypersurfaces in
$\C^2$ are self-dual. The non-quadratic CR-manifolds with
non-linearizable automorphisms are never self-dual.

If the complete solution of an ODE is known then the dual ODE can
be easily obtained by differentiating with respect to the
parameters and eliminating the variables $x,y$.

In the case of $y''=(y-xy')^3$ the complete solution is the
quartic
$$(y-c_1x)^2 - c_2^2x^2 - c_2=0$$

We find the dual ODE
\begin{equation}\label{dual}
y''=\frac{1-(y')^2}{x(y'+\sqrt{(y')^2-1})}.
\end{equation}

(Here we adopted $c_2$ as the new independent variable $x$ and
$c_1$ as the new dependent variable $y$.)

The family of solutions can be written in the form
$$(x-c_1)^2 - (y-c_2)^2=c_1^2.$$

The symmetries are generated by
$$\frac{\partial}{\partial y},\quad x\frac{\partial}{\partial x}+ y\frac{\partial}{\partial
y},\quad 2xy \frac{\partial}{\partial x}+ (x^2+y^2)
\frac{\partial}{\partial y}$$

The ODE (\ref{dual}) is equivalent to
$$\eta''+  \frac{2\eta'(1-\sqrt{\eta'})^2}{\xi-\eta}=0$$

from Lie's list of ODE with three symmetries. The equivalence is
established by $\xi=y+x$, $\eta=y-x$. In this notation the
infinitesimal automorphisms become
$$\frac{\partial}{\partial \xi}+ \frac{\partial}{\partial \eta} ,\quad
\xi\frac{\partial}{\partial \xi}+
\eta\frac{\partial}{\partial \eta},\quad \xi^2
\frac{\partial}{\partial \xi}+ \eta^2 \frac{\partial}{\partial
\eta}$$

The corresponding group is $\PSL(2,\C)$ acting by ``coupled''
M\"obius transformations on the complex $\xi$ and $\eta$ planes
$$(\xi,\eta) \mapsto \left(\frac{\alpha \xi +\beta}{\gamma \xi +\delta},
\frac{\alpha \eta +\beta}{\gamma \eta +\delta}\right).$$

\section{Shear invariant ODE with one non-isotropic symmetry}
\label{onenon} If a shear invariant ODE admits a non-isotropic
symmetry we may assume that, after a coordinate change $x=f(x^*,y^*)$, $y=g(x^*,y^*)$, it takes the form $\pd{x^*}{}$. Then the shear
becomes $\theta=\xi\pd{x^*}{} +\eta \pd{y^*}{}$ where $\xi=g \frac{\partial f^*}{\partial x}$, $\eta=g \frac{\partial g^*}{\partial x}$, and $(f^*,g^*)$ is the inverse coordinate change. We prove
\begin{lem}
If $\pd{x^*}{}$ and the shear $\theta$ are the only symmetries of the ODE then
$$[\pd{x^*}{}, \theta]=\mu \theta$$
\end{lem}

\begin{pf}
From
$$[\pd{x^*}{}, \theta]=\mu \theta + \nu \pd{x^*}{},$$
we conclude
\begin{align*}
\pd{x^*}{\xi}&=\nu + \mu \xi\\
\pd{x^*}{\eta}&= \mu \eta.
\end{align*}

Now, if $\mu=0$, then
\begin{align*}
\xi&=\nu x^* + K_1(y^*) \\
\eta&= K_2(y^*).
\end{align*}
We distinguish two subcases: If $K_2\equiv 0$ then $\pd{x}{g^*}\equiv 0$, and therefore $g=g(y^*)$ with $g(0)=0$. It follows that $\xi=0$ for $y^*=0$. Since $\theta$ vanishes exactly at one curve, this curve must coincide with $y^*=0$. Hence $\nu=0$.

In the second subcase $y^*=0$ is an isolated zero of $K_2$. Thus again $y^*$ is the only curve on which $\theta$ can vanish and therefore $\nu=0$.

Suppose now that $\mu\neq 0$. Then
\begin{align*}
\xi&=\frac{-\nu}{\mu} + K_1(y^*) \e^{\mu x^*}\\
\eta&= K_2(y^*) \e^{\mu x^*}.
\end{align*}
Again, either $K_2\equiv 0$ or $0$ is an isolated zero of $K_2$. Analogous arguments to the ones used above show that $\nu=0$ in this case as well.
\end{pf}

As in  Section \ref{oneiso} we consider the equation (\ref{inf}).
We conclude $\psi(y)=\alpha_0$ but now we assume that $\alpha_0
\not=0$. Then we can rescale the additional infinitesimal
automorphism is such a way that $\alpha_0=1$. Thus we look for an
infinitesimal automorphism of the form
$$(1+ (\phi(y)+a)x)\pd{x}{} + \phi(y)y \pd{y}{}.$$

From the component of degree 3 in $p$ and $1$ in $x$ we find
$$f_1=\frac{-\phi''}{2\alpha_0}= \frac{-\phi''}{2}.$$

The components of degree 3 in $p$ and $2$ resp. $3$ in $x$ yield
now the system
\begin{align}\label{f0}
-\frac{1}{6}(y\phi'''+3\phi'')\phi
-\frac{a}{6}\phi''&=f_0\\\nonumber (4\phi-y\phi'+2a)f_0+y\phi
f_0'&=0,
\end{align}

which is equivalent to the ODE
\begin{equation}\label{4ode}
 (y^2 \phi^{IV} + 8y\phi''' +12\phi'')\phi^2 +
a(3y\phi'''\phi+10 \phi''\phi - y\phi''\phi') +2a^2\phi''=0
\end{equation}

on $\phi$. We see immediately that $a= 0$ implies $\phi''=0$ and
therefore $f_0=f_1=0$. Assume $a\neq 0$.

The ODE (\ref{4ode}) yields the following equations on the
coefficients of an analytic solution $\phi(y)=\sum_{n=0}^\infty
\phi_n y^n$ :
\begin{align*}
\sum_{\beta=0}^{j-2} \sum_{\alpha=0}^{\beta}
\frac{(j-\beta+2)!}{(j-\beta-2)!} \phi_\alpha \phi_{\beta-\alpha}
\phi_{j+2-\beta} + 8 \sum_{\beta=0}^{j-1} \sum_{\alpha=0}^{\beta}
\frac{(j-\beta+2)!}{(j-\beta-1)!} \phi_\alpha \phi_{\beta-\alpha}
\phi_{j+2-\beta}+\\
+12\sum_{\beta=0}^{j} \sum_{\alpha=0}^{\beta}
\frac{(j-\beta+2)!}{(j-\beta)!} \phi_\alpha \phi_{\beta-\alpha}
\phi_{j+2-\beta} + 3a \sum_{\beta=0}^{j-1}
\frac{(j-\beta+2)!}{(j-\beta-1)!}  \phi_{\beta}
\phi_{j+2-\beta}+\\
+10a\sum_{\beta=0}^{j} \frac{(j-\beta+2)!}{(j-\beta)!}
\phi_{\beta} \phi_{j+2-\beta} -a \sum_{\beta=0}^{j}
\frac{(j-\beta+1)!}{(j-\beta-1)!} (\beta+1) \phi_{\beta+1}
\phi_{j+1-\beta}+\\
+2a^2 (j+2)(j+1) \phi_{j+2}=0
\end{align*}

It follows for $j\ge 0$
\begin{multline*}
(j+2)(j+1)(a+(j+3)\phi_0)(2a+(j+4)\phi_0) \phi_{j+2}+\\+
(j+2)(j+1)j (3a+2(j+3)\phi_0)\phi_1 \phi_{j+1}=\cdots
\end{multline*}
The dots indicate a sum, whose summands contain only factors
$\phi_n$ with $n\le j$ and at least one factor $\phi_n$ with $n\ge
2$.

Let $\phi_k$ with $k\ge 2$ be the first non-vanishing coefficient.
Then either
\begin{align*}
\phi_0&=-\frac{a}{k+1}, &\text{ or }\qquad
\phi_0&=-\frac{2a}{k+2},\\
\intertext{and, consequently,}
\phi_1&=\frac{a\phi_{k+1}}{(k-1)(k+1)\phi_k}, &\text{ or }\qquad
\phi_1&=\frac{2 a\phi_{k+1}}{(k-1)\phi_k},
\end{align*}

respectively.

For any parameter $a\neq 0$ related to the automorphism, we obtain
two series of solutions:

If $2a+(k+2)\phi_0=0$ (second option) then
$(a+(j+1)\phi_0)(2a+(j+2)\phi_0)\neq 0$ for all $j\ge k+2$ and
therefore all $\phi_j$ with $j\ge k+2$ can be obtained recursively
for given parameters $k, \phi_k\neq 0, \phi_{k+1}$.

If $a+(k+1)\phi_0=0$ (first option) then again all $\phi_j$ with
$j\ge k+2$ can be obtained recursively, except for $j=2k+2$. Here
an additional parameter $\phi_{2k+4}$ appears.

All other components in $\phi$ (and, thus, in $f_0,f_1$) can be
obtained recursively from
\begin{multline*}
(a+(j+k+3)\phi_0)(2a+(j+k+4)\phi_0) \phi_{j+k+2}+\\+
 (j+k)(3a+2(j+k+3)\phi_0)\phi_1 \phi_{j+k+1}+\\
+\sum_{\beta=0}^{j-k+2} \sum_{\alpha=0}^{\beta}
\frac{(j-\beta-k+4)!}{(j-\beta-k)!}
\frac{\phi_{k+\alpha}\phi_{k+\beta-\alpha}\phi_{j-\beta-k+2}}{(j+k+2)(j+k+1)}\\+a
\sum_{\alpha=0}^{j}
\frac{(3j-k-4\alpha+10)(j+2-\alpha)(j+1-\alpha)\phi_{k+\alpha}\phi_{j+2-\alpha}}{(j+k+2)(j+k+1)}=0
\end{multline*}

The convergence of the formal solutions can be proven by induction.

By applying a map of the form
$$x_1= \frac{c_1 x}{1-cy}, \qquad y_1=\frac{c_2 y}{1-cy}$$

we can renormalize a solution in such a way that $-a=\alpha_0=1$,
$f_{1,k-2}=-\frac{k(k-1)\phi_k}{2\alpha_0}=1$ and
$f_{1,k-1}=-\frac{k(k+1)\phi_{k+1}}{2\alpha_0}=0$. Thus, up to
equivalence, we obtain exactly two series of solutions, such that
a solution of the first series is determined by a non-negative
integer $k$ and a solution of the second series is determined by a
non-negative integer $k$ and a complex number $C$ which is related to $\phi_{2k+4}$. In all cases we
have
\begin{align*}
f_0(y)&=-\frac{1}{6}(y\phi'''+3\phi'')\phi +\frac{1}{6}\phi''\\
f_1(y)&=\frac{-\phi''}{2}
\end{align*}

with the additional symmetry
$$(1+ (\phi -1)x)\pd{x}{} + \phi y \pd{y}{},$$

where $\phi$ satisfies (\ref{4ode}) with initial conditions
$$\phi_0=\frac{1}{j+1},\quad \phi_1=0 \qquad \text{ or
}\qquad \phi_0=\frac{2}{j+2}, \quad \phi_1=0.$$

The first option corresponds to
$$y''=y^j(y-(x-c)y')^3,$$

which is obtained by shifting the ODE (\ref{vier1}) in the
$x$-direction by $c$. The parameter $c$ can be rescaled by
applying the additional isotropic automorphism.

The second option corresponds to shifts of (\ref{vier3})
$$y''=y^jy'(y-(x-c)y')^2 + C y^{2j+2}(y-(x-c)y')^3.$$

In the special case $C=0$ we deduce $f_0\equiv 0$ and
$$(y\phi'''+3\phi'')\phi -\phi''=0.$$

In terms of $f_1$ the latter equation becomes
$$\left( \frac{f_1}{3f_1+yf_1'}\right)''=-2f_1.$$

\section{Shear invariant ODE with two additional
symmetries}\label{twosy}

If a shear invariant ODE has two additional symmetries then either
one of them can be chosen to be isotropic or both give rise to a
transitive sub-semigroup on $\C^2$. According to the results of
Section \ref{oneiso} the first case leads to three particular
series of ODE, which have only isotropic symmetries. The only ODE
(up to equivalence) with two additional isotropic symmetries is
$$y''=(y-xy')^3.$$
Consider the second case. We may assume that there is an
infinitesimal non-isotropic automorphism $\sigma$ in the direction
of the line of fixed points of the shear $\theta$. Without loss of
generality we have then
$$\sigma=\pd{x}{}, \qquad \theta=(y+ a)\pd{x}{} + b \pd{y}{},$$

where $a(x,y), b(x,y)$ are of at least second order. But then
$$[\sigma,\theta]=\lambda \theta$$
because $\theta$ is the only isotropic symmetry. According to the
results of Section \ref{onenon} we conclude that the ODE must be a
shift of (\ref{vier1}) or (\ref{vier3}). Again, only
$$y''=(y-(x-c)y')^3$$

has three-dimensional symmetry.

\section{Non-linearizable automorphisms of elliptic
CR-manifolds}\label{nonlin}

In \cite{ES03} we proved that the phenomenon of non-linearizable
isotropy takes place on a whole complex curve. As a consequence of
the classification results from above we prove here the following
converse statement for an elliptic CR manifold $M$ with the
additional property that all infinitesimal automorphisms are
globally defined.

\begin{theo}
Let $M$ be an elliptic CR-manifold with non-linearizable isotropy
group at $p\in M$. If $M$ is neither equivalent to the quadric
$$w_1-\bar{w}_2-z_1\bar{z}_2=0$$
nor to the quartic
$$w_1 +w_1^2\bar{z}_2^2-(\bar{w}_2-z_1\bar{z}_2)^2 =0 $$
then there exists a neighbourhood $U$ of $0$ such that $\Aut_q M$
is linearizable for $q\in U$ outside a complex curve $\gamma$.
\end{theo}

\begin{pf}
Let $q\in M$ be a point with non-linearizable isotropy. If $M$ is
not equivalent neither to the quadric nor to the quartic then
there is a single shear at $q$, which either coincides with the
single shear in $0$ or it provides an additional symmetry at $0$.
In the first case $q$ is a fixed point of $y\pd{x}{}$ and
therefore belongs to $\gamma=\{y=p=0\}$.

In the second case $M$ corresponds to one of the ODEs listed
above. But then only the shear has non-isolated fixed points
outside $0$. All these fixed points belong to $\{y=p=0\}$.
\end{pf}

The quartic can be characterized by the following property.

\begin{prop}
The set of points at the quartic with non-linearizable isotropy is
the complex hypersurface $\Gamma=\{y=xp\}$.
\end{prop}

\begin{pf}
The mapping
$$(x, y) = (a (x_1+1) +c y_1, b(x_1+1) +dy_1)$$
takes the the point $(a,b, \frac{b}{a})$ to $(0,0,0)$ and the ODE
$y''=(y-xy')^3$ again to an ODE that admits a shear, namely to
$$y''=(ad-bc)^2(y-(x+1)y')^3.$$
Since the orbit of $0$ under these mappings is the hypersurface
$\Gamma$, at all points of $\Gamma$ the isotropy group is
non-linearizable.

We show that the isotropy of the quartic is linearizable (even
trial) at any point outside $\Gamma$. Any infinitesimal
automorphism of the quartic at $0$ has the form
$$(\alpha x +\beta y) \pd{x}{} + (\delta x -\alpha y)\pd{y}{} + (\delta - 2\alpha p
-\beta p^2)\pd{p}{}.$$

If the discriminant $\Delta=\alpha^2+\beta\delta$ is different
from $0$ then fixed points occur only for $x=y=0$. If the
discriminant vanishes we distinguish the two subcases $\beta\neq
0$ and $\beta=0$. In the first subcase we find fixed points for
$\alpha x +\beta y=0$, $p=-\frac{\alpha}{\beta}$. This implies
$y-xp=0$. If $\beta=0$ we conclude $\alpha=\delta=0$. Then only
the identical automorphism has fixed points other than $0$.
\end{pf}

\end{document}